     \def\({\left(}                  
     \def\){\right)}      \def\e{\varepsilon}    
     \def\[{\left[}       
     \def\]{\right]}

     \def\<{\langle}                 
     \def\>{\rangle}                 \def\wt{\widetilde}
                 \def\sbs{\subset}

\def\dim{\operatorname{dim}}

                     \def\sbs{\subset}

\def\span{\operatorname{span}}

\documentclass[12pt,oneside]{amsproc}
\usepackage[T2A]{fontenc}
\usepackage[cp1251]{inputenc}
\usepackage{amssymb}
\usepackage{amsthm}

\title{On operators with bounded approximation property}
\author{Oleg Reinov}
\address{ St. Petersburg State University,
Saint Petersburg, RUSSIA.}
\email{orein51@mail.ru}

\thanks{
AMS Subject Classification 2010: 46B28 Spaces of operators; tensor products; approximation properties.
}
\thanks{${ }$ Key words: 46B28  Approximation of operators; bounded approximation property.}

\begin{document}

                          $$ {} $$
\maketitle


\begin{abstract}
It is known that any separable Banach space with BAP is a complemented subspace 
of a Banach space with a basis. We show that every operator with bounded approximation property, 
acting from a separable Banach space,   can be factored through a Banach space with a basis.
 \end{abstract}

\vskip 0.3cm

\centerline{\bf \S1. Lemmas}

 \vskip 0.23cm

{\bf Definition 1.1.}\,
Let $T\in L(X,W),$  $C\ge1.$ We say that $T$ has the C-BAP (C-bounded approximation property)
if for every compact subset $K$ of $X,$ for any $\e>0$ there exists a finite rank operator $R: X\to W$
such that $||R||\le C\, ||T||$ and $\sup_{x\in K} ||Rx-Tx||\le \e.$ The operator $T$ has the BAP if
it has the C-BAP for some $C\in[1,\infty).$

\medskip

{\bf Lemma 1.1.}\,
$T$ has C-BAP iff for any finite family $(x_k)_{k=1}^N\sbs X,$ for any $\e>0$
there exists a finite rank operator $R: X\to W$ such that
$||R||\le C||T||$ and $\sup_{1\le k\le N} ||Rx_k-Tx_k||\le \e.$

{\it Proof}.\,
We may (and do) assume that $||T||=1.$
Fix a compact subset $K\sbs X$ and $\e>0.$ Let $\e_0:= \e/(2+C),$ $(x_k)_{k=1}^M$
be an $\e_0$-net for $K$ in $X,$ $R\in X^*\otimes W,$ $||R||\le C$ and $\sup_{1\le k\le M} ||Rx_k-Tx_k||\le \e_0.$
Take an $x\in K,$ and let $x_k$ be such that $||x-x_k||\le \e_0.$ Then
$||Tx-Rx||\le ||x-x_k||+ ||Tx_k-Rx_k|| +||Rx_k-Rx||\le \e_0+\e_0+C\e_0=\e.$
\smallskip

{\bf Lemma 1.2.}\,
Let $X, W$ be Banach spaces, $X$ being separable, and $T\in L(X, W).$
$T$ has C-BAP iff  there exists a sequence $(Q_l)_{l=1}^\infty$ of finite rank operators
from $X$ to $W$ such that

$1)$\, for every $x\in X$ the series $\sum_{l=1}^\infty Q_lx$ converges and
$$
 Tx=\sum_{l=1}^\infty Q_lx,\ x\in X;
$$

$2)$\, $\sup_N ||\sum_{l=1}^N Q_l||\le C\, ||T||.$
\smallskip

{\it Proof}.\,
Since $X$ is separable, there exists a sequence $(x_k)_1^\infty$ which is dense
in the closed unit ball $\bar B_1(0)$ of $X.$
Suppose  as above that $||T||=1$ and $T$ has the C-BAP, that is
for any finite set $F\sbs X,$ for every $\e>0$ there is a finite rank operator $R:X\to W$
such that $||R||\le C$ and $\sup_{f\in F} ||Rf-Tf||\le \e.$
Put, for $N=1,2, \dots$, $F_N:= \span (x_k)_{k=1}^N;$  $F_N\sbs F_{N+1} \dots.$
For each $N,$ let $R_N$ be a finite rank operator from $X$ to $W$ with the properties that

(i)\, $||R_N||\le C$ and

(ii)\, $\sup_{1\le n\le N} ||R_Nx_n-Tx_n||\le 1/2^{N+1}.$

If $n\in \Bbb N$ then for every $N\ge n$ one has
$$
 (iii)\ \ ||R_Nx_n-Tx_n||\le \frac1{2^{N+1}}
$$
and, therefore, for a fixed $x_n$
$$
 R_Nx_n \rightarrow Tx_n
$$
as $N$ tends to $\infty.$

Now, fix $\e>0$ and let $\delta>0$ be such that $C\delta+\delta<\e.$
For $x\in \bar B_1(0),$  take an $x_n$ with $||x_n-x||<\delta.$
Then there is an $N_0$ so that for $N\ge N_0$
$$
 ||R_Nx-Tx||\le ||R_N||\, ||x_n-x|| +||R_Nx_n-Tx_n||+ ||Tx_n-Tx||\le C\delta +||T||\delta<\e.
$$
Thus, if $x\in X$ then $R_Nx \to Tx$ as $N\to +\infty.$

To finish the proof of the "only if" part, we apply
\smallskip

{\bf Lemma 1.3.}\,
Let $X, W$ be any Banach spaces,  $C\ge1$ and $T\in L(X, W).$
Suppose that

(*)\, there exists a sequence $(S_N)_{N=1}^\infty$ of finite rank
operators from $X$ to $W$ such that
if $x\in X$ then $S_Nx \to Tx$ as $N\to +\infty$ and $||S_N||\le C||T||$ for every $N.$

Then
 there exists a sequence $(Q_l)_{l=1}^\infty$ of finite rank operators
from $X$ to $W$ such that

$1)$\, for every $x\in X$ the series $\sum_{l=1}^\infty Q_lx$ converges and
$$
 Tx=\sum_{l=1}^\infty Q_lx,\ x\in X;
$$

$2)$\, $\sup_N ||\sum_{l=1}^N Q_l||\le C\, ||T||.$
\smallskip

{\it Proof}.\,
We assume again that $||T||=1.$
Put $Q_1:=S_1, Q_l:=S_l-S_{l-1}$ for $l>1,$ so that
$$
  S_N= S_1+(S_2-S_1)+\dots +(S_{N-1}-S_{N-2})+ (S_N-S_{N-1})=
  Q_1+Q_2+\dots+Q_N.
$$
It follows that

$$(1)\ Tx=\sum_{l=1}^\infty Q_lx\ \ \forall\, x\in X$$
and
$$(2)\ \sup_{N}||\sum_{l=1}^NQ_j||=\sup_N ||S_N||\le C.$$

The "if" part of the proof of Lemma 1.2 follows from
\smallskip

{\bf Lemma 1.4.}\,
Let $X, W$ be any Banach spaces, $C\ge1$ and $T\in L(X,W).$
If there exists a sequence $(R_N)_1^\infty$
of finite rank operators from $X$ into $W$ which converges pointwise to $T$
and such that $||R_N||\le C||T||$ for all $N$
then $T$ has the C-BAP.
\smallskip

{\it Proof}.\,
Indeed, let $R_Nx\to Tx$ for every $x\in X$ (and $||R_N||\le C||T||$). Fix $\e>0$
and a compact subset $K\sbs X.$ Put $\e_0:=\e (||T||+1+C||T||)^{-1}.$
Take a finite $\e_0$-net $F\sbs X$ for $K$ and consider $R_{N_0}$ such that
$\sup_{f\in F} ||R_{N_0}f- Tf||\le \e_0.$ Then, for any $x\in K$ there is
an $f_0\in F$ with $||f_0-x||\le\e_0,$ and one has :
$$
  ||Tx-R_{N_0}x||\le ||T||\, \e_0+\e_0+||R_{N_0}||\e_0\le
  \e_0 (||T||+1+C||T||)=\e.
$$
\smallskip

{\bf Corollary 1.1.}\,
If $X$ is separable and $T\in L(X,W),$ then $T$
has the C-BAP iff there exists a sequence $(R_N)_1^\infty$
of finite rank operators from $X$ into $W$ which converges pointwise to $T$
and such that $||R_N||\le C||T||$ for all $N.$
\smallskip

{\bf Corollary 1.2.}\,
If $X$ is separable and $T\in L(X,W),$ then $T$
has the BAP iff there exists a sequence
of finite rank operators from $X$ to $W$  convergent  to $T$ pointwise.

\bigskip

\centerline{\bf \S2. Theorem.}

\bigskip

Now, we redenote some objects from \S1. Let $X, W$ be any Banach spaces,
$T\in L(X,W)$ and $T$
 possesses the property (*) from Lemma 1.3.
 Consider the sequence $(Q_l)_{l=1}^\infty,$ given by assertion of Lemma 1.3,
and put $A_p:=Q_p$\, $(p=1,2,\dots)$ and $K:=C (\ge1),$ assuming that $||T||=1.$
We are now in notations (partially) of the paper [1].          
\smallskip

{\bf Theorem 2.1.}\,
If $T:X\to W$ has the property (*), then there exist a Banach space $Y$
with a Schauder  basis and two operators $\wt A: X\to Y$ and $j: Y\to W$  
so that $T=j\wt A.$
\smallskip

{\it Proof}.\,
In the above notation (assuming that $||T||=1$), we have:
$$
 Tx=\sum_{p=1}^\infty A_px,\ \forall\, x\in X;\  \ A_p\in X^*\otimes W,\  \sup_{n\in\Bbb N} ||\sum_{p=1}^n A_p||\le K
$$
(note that for every $n$ \, $||A_n||\le ||\sum_{p=1}^n A_p- \sum_{p=1}^{n-1} A_p||\le 2K$).
Let $E_p=A_p(X)\sbs W,$ \, $m_p:=\dim E_p$ for $p\ge1$ and $m_0=0.$
We will proceed as in [1].                     

By Auerbach, there exist one-dimensional operators $B_j^{(p)}: E_p\to E_p$
with $||B_j^{(p)}||=1$ for $j=1,2,\dots, m_p,$
and so that
$$
 \sum_{j=1}^{m_p} B_j^{(p)}(e)=e,\ e\in E_p.
$$
Set $C_i^{(p)}:= \dfrac1{m_p}\, B_j^{(p)}$ for $i=rm_p+j$ (where $r=0,1,\dots, m_p-1; j=1,2,\dots, m_p).$
Then, for $e\in E_p.$
$$
  \sum_{i=1}^{m_p^2} C_i^{(p)} (e)=m_p \cdot \sum_{j=1}^{m_p} \frac1{m_p} B_j^{(p)}e=
  \sum_{j=1}^{m_p}  B_j^{(p)}e=e.
$$
 Also, for any $q\ge1, q\le m_p^2$ and some $l< m_p$ and $k\le m_p$ we have:
 $$
   ||\sum_{i=1}^q C_i^{(p)}||= ||\sum_{i=1}^{lm_p} C_i^{(p)} +\sum_{lm_p+1}^{lm_p+k} C_i^{(p)}|| \le
       l\cdot \frac1{m_p} \, ||\sum_{j=1}^{m_p} B_j^{(p)}|| + \frac1{m_p}\cdot ||\sum_{j=1}^{k} B_j^{(p)}|| \le 1+1=2.
 $$

 Now, let
 $$\wt A_s:= C_i^{(p)}A_p$$
  for $p\in \Bbb N,$ $i=1,2,\dots, m_p^2$ and $s=m_0^2+m_1^2+\dots +m_{p-1}^2+i.$
  1-dimensional operator $\wt A_s$ maps $X$ into $E_p\sbs W$ in the following way:
  $$
   \wt A_s: X\stackrel{A_p}{\rightarrow} E_p=A_p(X)\stackrel{C_i^{(p)}}{\rightarrow} E_p (\sbs W).
$$
 Since, for any $n\in\Bbb N,$ for some $k$ and $r\le m_k^2$

 $$
  \sum_{s=1}^n \wt A_s=\sum_{p=1}^{k-1} \sum_{i=1}^{m_p^2} C_i^{(p)}A_p + \sum_{i=1}^{r} C_i^{(k)}A_k,
 $$
 we get that
 $$
  || \sum_{s=1}^n \wt A_s||\le ||\sum_{p=1}^{k-1} A_p|| + ||\sum_{i=1}^{r} C_i^{(k)}A_k||\le K+2||A_k||\le 5K.
 $$
 (To get an estimation "$4K$" as in [1], it it enough to consider simultaniousely, in the center,
 the given sum and the sum like $||\sum_{p=1}^{k} A_p|| + ||\sum_{i=r}^{m_k^2} C_i^{(k)}A_k||$).
 
 Since, for every $x\in X,$ $A_kx\to 0$ as $k\to\infty,$ we have:
 $$
  \lim_{n\to\infty}  \sum_{s=1}^n \wt A_s x= \lim_{N\to\infty} \sum_{p=1}^{N} A_px =Tx.
 $$

Now, consider the space 
$$Y:=\{(y(s))_{s=1}^\infty:\  y(s)\in \wt A_s(X), \, \sum_{s=1}^\infty y(s)\ \mathrm{converges\ in }\ W\}.$$
Set $|||(y(s))_1^\infty|||:= \sup_n ||\sum_{s=1}^n y(s)||_W.$
Note that $(\wt A_s(x))_{s-1}^\infty \in Y$ for every $x\in X,$ $\sum_{s-1}^\infty \wt A_s(x)= Tx$ and
$|||(\wt A_s(x))_{s=1}^\infty|||_Y\le 5K||x||_X.$
Therefore, the map $\wt A:X\to Y,$ defined by $\wt A(x)=(\wt A_s(s))_{s=1}^\infty,$
is linear and continuous (and $||\wt A||\le 5K$).
Let $j: Y\to W$ be the natural map which takes $(y(s))_{s=1}^\infty$ to $\sum_{s=1}^\infty y(s).$
Since 
$$||\sum_{s=1}^\infty y(s)||_W = \lim_{N} ||\sum_{s=1}^N y(s)||_W\le \sup_n ||\sum_{s=1}^n y(s)||_W,$$
then $||j||_{L(Y,W)}\le 1.$
Therefore, $Tx=j\wt A: X\to Y\to W.$
It remains now to consider the space $Y.$

For each $s,$ let $\wt y_s\in \wt A_s$ be of norm 1. If $(y_s)_{s=1}^\infty \in Y,$
then $y_s=c_s\wt y_s.$ 
Define $\bar y_s\in Y$ by $\bar y_s(t)=0$ for $t\neq s$ and $\bar y_s(s)=\wt y_s$ \, $(s=1,2,\dots).$
Then,  every $y=(y_s)\in Y$ is of type $\sum_{s=1}^\infty c_s\bar y_s,$
if we consider $(\bar y_s)_{s=1}^\infty$ as a basis in Y. And this basis is monotone:
for all  scalars $(c_s)$ we have that
$$
 |||\sum_{s=1}^m c_s \bar y_s|||\le |||\sum_{s=1}^{m+1} c_s \bar y_s|||
$$
(by definition of the norm in $Y).$
Finally, the space $Y$ is Banach (cf. [2,p. 18, Prop. 3.1]).
\medskip

{\it Remark}.\,
The theorem just obtained is a spade-theorem for some futher investigations in a next paper.

\medskip

{\bf Corollary 2.1.}\,
If $X$ is separable and $T\in L(X,W)$ then $T$ has the bounded approximation property if and only if
$T$ can be factored through a Banach space with a basis.





\end{document}